\documentclass[letterpaper, 10 pt, conference]{ieeeconf}  %
\IEEEoverridecommandlockouts                              %
\usepackage[utf8]{inputenc} %
\usepackage[T1]{fontenc}    %
\usepackage{balance}
\usepackage[english]{babel}
\usepackage{lmodern, verbatim}
\usepackage{calc}
\usepackage{hyperref}       %
\usepackage{url}            %
\usepackage{booktabs}       %
\usepackage{amsfonts}       %
\usepackage{nicefrac}       %
\usepackage{amssymb}
\usepackage{amsmath}
\usepackage[pdftex]{graphicx}

\usepackage{subfig}
\usepackage{color}%

\usepackage{amsthm}

\usepackage{cleveref}

\newcommand{\RR}{\mathbb{R}}
\newcommand{\R}{\RR} %
\renewcommand{\b}{\mathbf}
\DeclareMathOperator*{\Id}{Id}
\newcommand{\cH}{{\mathcal H}}
\newcommand{\cL}{{\mathcal L}}
\newcommand{\HK}{{\mathcal H}_K}

\providecommand{\noridx}[2]{\lVert{#1}\rVert_{#2}}

\providecommand{\braidx}[3]{\left\langle{#1},{#2}\right\rangle_{#3}}
\providecommand{\scalidx}[3]{\left({#1},{#2}\right)_{#3}}

\newcommand{\1}{\mbox{1\hspace{-1mm}I}}

\theoremstyle{plain}
\newtheorem{thm}{Theorem}

\newtheorem{claim}{Claim}
\newtheorem{prop}{Proposition}%

\title{\LARGE \bf Operator-valued Kernels and Control of Infinite dimensional Dynamic Systems
}

\author{Pierre-Cyril Aubin-Frankowski and Alain Bensoussan%
	\thanks{This work was supported by the National Science Foundation under grant NSF-DMS-1905449 and grant from the SAR Hong Kong RGC GRF 14301321, and by the European Research Council (grant REAL 947908).}%
	\thanks{P.-C.\ Aubin-Frankowski is with INRIA and Département d’Informatique, École Normale Supérieure, PSL Research University, Paris, France.
		{\tt\small pierre-cyril.aubin@inria.fr}}%
	\thanks{A. Bensoussan is with the International Center for Decision and Risk Analysis, Jindal School of Management and University of Texas at Dallas, also with the School of Data Science, City University Hong Kong.
		{\tt\small axb046100@utdallas.edu}}%
}

\begin{document}

	\maketitle
	\thispagestyle{empty}
	\pagestyle{empty}

	\begin{abstract}
		The Linear Quadratic Regulator (LQR), which is arguably the most classical problem in control theory, was recently related to kernel methods in \cite{aubin2020hard_control} for finite dimensional systems. We show that this result extends to infinite dimensional systems, i.e.\ control of linear partial differential equations. The quadratic objective paired with the linear dynamics encode the relevant kernel, defining a Hilbert space of controlled trajectories, for which we obtain a concise formula based on the solution of the differential Riccati equation. This paves the way to applying representer theorems from kernel methods to solve infinite dimensional optimal control problems.
	\end{abstract}

	\section{Introduction}
	We follow the paper \cite{aubin2020hard_control} of one of the authors in which a connection was made between the LQR problem for finite dimensional systems and the important theory of Reproducing Kernel Hilbert Spaces (RKHSs), introduced by N.\ Aronszajn \cite{Aronszajn1943,aronszajn50theory}. This theory, under the name of kernel methods, has been successfully used in nonparametric statistics \cite{berlinet04reproducing}, and more recently raised a strong interest in the machine learning community \cite{steinwart08support}.  We refer to \cite{aubin2022Kalman} for connections between kernels as covariances and estimation problems. This paper extends the results of \cite{aubin2020hard_control} in two directions. First we obtain new explicit formulas for the kernel associated with the LQR and apply them to recover the solution of the unconstrained LQR through kernel principles, and second we tackle infinite dimensional linear systems. In doing this, we need to work with operator-valued kernels, which have been mostly used in a different context, namely supervised learning for functional data analysis (see \cite{Micchelli2005,caponetto2008universal,kadri2016} and references therein). Control theory for infinite dimensional dynamical systems, with application to the control of distributed parameter systems has been strongly developed by J.\ L.\ Lions \cite{lions1971optimal,lions1988controlabilite}. We use his notation, and many of his results.%

	The paper is organized as follows. In \Cref{sec:LQR}, we set the context, recall the solution of the LQR and its connection with the differential Riccati equation. We identify the related operator-valued kernel in \Cref{theo3-1} presented in \Cref{sec:kernel_LQR}, wherein the role of the space of trajectories with null initial condition is then discussed, and exemplified on the heat equation with distributed control. In \Cref{sec:optimal_control}, we apply the kernel to formally solve optimal control problems, recovering the LQR solution as a side result.
	
	\section{Linear Quadratic Regulator for infinite dimensional dynamic systems}\label{sec:LQR}
	
	\subsection{Notation and statement of the problem}
	Let $(V,\noridx{\cdot}{V})$ and $(H,\noridx{\cdot}{H})$ be two separable Hilbert spaces. We assume that
	$V\subset H$,  with continuous injection. Identifying $H$ to its
	dual, we have also the inclusion $H\subset V'$  with continuous injection, where $V'$ is the dual of $V$. We then consider a family
	$A:t\mapsto A(t)$ of continuous linear operators $A(t)\in\cL(V,V')$
	and we assume a classical coercivity condition (e.g.\ \cite{lions1971optimal} pp.100 and 272)
	
	\vspace{-8mm}
	\begin{gather}
		\exists \alpha>0,\beta\in\R,\, \forall z\in V,\nonumber\\ \braidx{A(t)z}{z}{V'\times V}+\beta \noridx{z}{H}^{2}\geq\alpha\noridx{z}{V}^{2}, \text{a.e. }t\in [t_0,T].\label{eq:2-1}
	\end{gather}%
	Let $U$ be another Hilbert space, and take $B(\cdot)\in L^{\infty}(t_0,T;\cL(U,H)).$ 
	If $u(\cdot)\in L^{2}(t_0,T;U)$ is a control, the state evolution in $H$
	of a linear system with initial condition $y_{0}\in H$ is described by the equation 
	\begin{equation}
		\tfrac{dy}{dt}+A(t)y(t)=B(t)u(t), \, y(t_0)=y_{0}\in H\label{eq:2-2}
	\end{equation}
	The solution $y(\cdot)$ of the linear system (\ref{eq:2-2}) belongs
	to $L^{2}(t_0,T;V)$, whereas $ \nicefrac{dy}{dt}(\cdot)\in L^{2}(t_0,T;V').$ The Linear
	Quadratic Regulator (LQR) is the problem of minimizing 
	
	\vspace{-3mm}
	\begin{equation}
		J(u(\cdot))=\int_{t_0}^{T} [\scalidx{M(t)y(t)}{y(t)}{H}+ \scalidx{N(t)u(t)}{u(t)}{U}]dt,\label{eq:2-3}
	\end{equation}
	where $M(\cdot)\in L^{\infty}(t_0,T;\cL(H,H))$, $N(\cdot)\in L^{\infty}(t_0,T;\cL(U,U))$
	and $M(t),N(t)$ are self-adjoint with 
	
	\vspace{-5mm}
	\begin{equation}
		M(t)\geq0,\quad N(t)\geq\nu \Id\nolimits_U, \; \text{ for some }\nu>0.\label{eq:2-4}
	\end{equation}
	For PDE control over an open bounded set of $\R^d$ with regular boundary, the assumption of boundedness of $B$ allows us to consider controls in the domain, but not control on the boundary (\cite{Bensoussan2007} Part IV-1), see for instance \cite{Pritchard1987} for the unbounded LQR. On the other hand, the coercivity assumption \eqref{eq:2-1} over $A(t)$ is satisfied by many parabolic or hyperbolic equations, such as the heat or the wave equations. We refer to (\cite{lions1971optimal} Section 3.4, p23) and (\cite{Bensoussan2007} Part IV-1:7.2, p415) for more examples. We only focus here on parabolic equations in first-order form \eqref{eq:2-2} and leave second-order hyperbolic or unbounded $B(t)$ settings to further research which would require a mild adaptation.
	
	\subsection{Solution of the linear quadratic regulator}
	
	The functional $u(\cdot)\mapsto J(u(\cdot))$ is quadratic and strictly convex.
	It has a unique minimum $u(\cdot)$, which is computed as follows: the
	forward-backward system of equations 
	
	\vspace{-5mm}
	\begin{gather}
		\tfrac{dy}{dt}+A(t)y(t)+B(t)N^{-1}(t)B^{*}(t)p(t)=0,\label{eq:2-5}\\
		-\tfrac{dp}{dt}+A^{*}(t)p(t)-M(t)y(t)=0, \nonumber \\
		y(t_0)=y_{0},\:p(T)=0, \nonumber
	\end{gather}
	has a unique solution $y(\cdot)\in L^{2}(t_0,T;V)$, with $\tfrac{dy}{dt}\in L^{2}(t_0,T;V')$, $p(\cdot)\in L^{2}(t_0,T;V)$ and $\tfrac{dp}{dt}\in L^{2}(t_0,T;V').$
	Moreover, we have the decoupling property 
	\begin{equation}
		p(t)=P(t)y(t)\label{eq:2-6}
	\end{equation}
	in which $P(t)$ $\in\cL(H;H)$ is symmetric and positive semidefinite. The
	operator $P(t)$ is defined by solving a system similar to (\ref{eq:2-5}) for each $t\in[t_0,T]$ and $h\in H$
	\begin{gather}
		\tfrac{d\xi}{ds}+A(s)\xi(s)+B(s)N^{-1}(s)B^{*}(s)\eta(s)=0,\label{eq:2-7}\\
		-\tfrac{d\eta}{ds}+A^{*}(s)\eta(s)-M(s)\xi(s)=0,\; \forall s\in(t,T), \nonumber \\
		\xi(t)=h,\:\eta(T)=0, \nonumber
	\end{gather}
	and then setting
		
	\vspace{-7mm}
	\begin{equation}
		\eta(t)=P(t)h.\label{eq:2-8}
	\end{equation}
	The operator $P(t)$ has the following regularity property: 
	\begin{claim}\label{eq:2-80}
		If $\varphi(\cdot)\in L^{2}(t_0,T;H)$ satisfies $\tfrac{d\varphi}{dt}+A(t)\varphi(t)\in L^{2}(t_0,T;H)$, then $\Psi(t)=P(t)\varphi(t)$ satisfies $-\tfrac{d\Psi}{dt}+A^{*}(t)\Psi(t)\in L^{2}(t_0,T;H)$, and 
		\begin{multline}
			-\tfrac{d\Psi}{dt}+A^{*}(t)\Psi(t)+P(t)\big[\tfrac{d\varphi}{dt}+A(t)\varphi(t)\\+B(t)N^{-1}(t)B^{*}(t)\Psi(t)\big]=M(t)\varphi(t).
		\end{multline}
	\end{claim}
	This formally can be written as 
	
	\vspace{-7mm}
	\begin{multline}
		-\tfrac{dP}{dt}+P(t)A(t)+A^{*}(t)P(t)\\ +P(t)B(t)N^{-1}(t)B^{*}(t)P(t)=M(t), \; P(T)=0.\label{eq:2-82}
	\end{multline}
	The optimal state $y(\cdot)$ for the LQR control problem (\ref{eq:2-2}),
	(\ref{eq:2-3}) is solution of the equation
	\begin{equation}
		\tfrac{dy}{dt}+(A(t)+B(t)N^{-1}(t)B^{*}(t)P(t))y(t)=0, \; y(t_0)=y_{0}.\label{eq:2-9}
	\end{equation}
	and the optimal control $u(\cdot)$ is given by 
	\begin{equation}
		u(t)=-N^{-1}(t)B^{*}(t)P(t)y(t).\label{eq:2-10}
	\end{equation}
	We will use in the sequel the semi-group (a.k.a.\ evolution family) associated with the operator
	$A(t)+B(t)N^{-1}(t)B^{*}(t)P(t)$,  denoted by $\Phi_{A,P}(t,s)$,   for $t>s$, i.e. $\Phi_{A,P}(s,s)=\Id\nolimits_H$ and
	
	\vspace{-5mm}
	\begin{multline}\label{eq:def_PhiAP}
		\partial_t\Phi_{A,P}(t,s)+\\(A(t)+B(t)N^{-1}(t)B^{*}(t)P(t))\Phi_{A,P}(t,s)=0.
	\end{multline}
	Note that if $M(\cdot)\equiv 0$, then $P(\cdot)\equiv 0$ and $\Phi_{A,P}(t,s)=\Phi_{A}(t,s)$, the semi-group associated with the operator $A(t)$, i.e.\ $\partial_t \Phi_{A}(t,s)+A(t)\Phi_{A}(t,s)=0$ and $\Phi_{A}(t,t)=\Id_H$.
	
	\section{Kernel associated with the Linear Quadratic Regulator} \label{sec:kernel_LQR}
	
	\subsection{$\cL(H;H)$-valued reproducing kernel}
	
	Owing to the results of \cite{Senkene1973} and Theorem 2.12 in \cite{burbea1984banach}, given a Hilbert space $(\cH,\noridx{\cdot}{\cH})$ of functions of $[t_0,T]$ into a separable Hilbert space $H$, such that the strong topology of $\cH$ is stronger than pointwise convergence, there exists a unique function $K:s,t\in[t_0,T]\mapsto K(s,t)\in\cL(H;H)$, which is called the reproducing $\cL(H;H)$-valued kernel of $\cH$, satisfying
	\begin{align}
		K(\cdot,t)z\in\cH,\; \forall\,t\in[t_0,T],\,z\in H, \label{eq:3-100}\\
		\scalidx{y(t)}{z}{H}=\scalidx{y(\cdot)}{K(\cdot,t)z}{\cH},\label{eq:rep-prop}\\ \forall\, y(\cdot)\in \cH,\,t\in[t_0,T],\,z\in H. \nonumber
	\end{align}
	Conversely a kernel $K$ satisfying \eqref{eq:3-100}-\eqref{eq:rep-prop} characterizes $\cH$. We then say that $\cH$ is a $H$-valued Reproducing Kernel Hilbert space (RKHS), and that $K(\cdot,\cdot)$ is the $\cL(H;H)$-valued kernel associated with $\cH$. 
	
	Following the initial work of one of the authors \cite{aubin2020hard_control}, we want to show that the LQR problem \eqref{eq:2-2}-\eqref{eq:2-3} can be associated with an $\cL(H;H)$-valued kernel. We use then the notation $\HK$ to emphasize the role played by the kernel.
	
	\subsection{Hilbert space of controlled trajectories}
	
	We consider the subset $\cH$ of $L^{2}(t_0,T;H)$ defined as
	follows
	\begin{multline}
		\cH=\{y(\cdot)\in L^{2}(t_0,T;H)\, | \,\tfrac{dy}{dt}+A(t)y(t)=B(t)u(t),\\ \text{with}\;u(\cdot)\in L^{2}(t_0,T;U)\}.\label{eq:3-1}
	\end{multline}
	There is not necessarily a unique choice of $u(\cdot)$ for a given $y(\cdot) \in \cH$ (for instance if $B(t)$ is not injective for some $t$). Therefore, with each $y(\cdot) \in \cH$, we associate the control $u(\cdot)$ having minimal norm based on the pseudoinverse of $B(t)^{\ominus}$ of $B(t)$ for the $U$-norm $\|\cdot\|_{N(t)}:=\|N(t)^{1/2}\cdot\|_U$:
	\begin{align}
		u(t) = B(t)^{\ominus}[\tfrac{dy}{dt} + A(t) y(t)] \, \text{ a.e.\ in} \,[t_0,T], \label{def_u_as_X-X'} %
	\end{align}
	whence $u(\cdot)$ minimizes $\int_{t_0}^{T} \scalidx{N(t)u(t)}{u(t)}{U}dt$ among the controls admissible for $y(\cdot)\in\cH$. We consequently equip $\cH$ with the norm 
	\begin{align}
		&\noridx{y(\cdot)}{\cH}^{2}=\scalidx{y(t_0)}{J_0 y(t_0)}{H}\nonumber\\
		&\hspace{2mm}+\int_{t_0}^{T} [\scalidx{M(t)y(t)}{y(t)}{H}+ \scalidx{N(t)u(t)}{u(t)}{U}]dt,\label{eq:3-2}
	\end{align}
	for some positive semidefinite $J_0\in  \cL(H;H)$ such that
	\begin{equation}
		J_0+P(t_0) \text{ is invertible.}\label{eq:invertibility}
	\end{equation}
	where $P(\cdot)$ is defined as in \eqref{eq:2-82}. Then $\cH$ has the structure of a Hilbert space. The extra term $J_0$ in \eqref{eq:3-2}, as compared to \eqref{eq:2-3}, is required to define a norm.\footnote{The condition $\noridx{y(\cdot)}{\cH}^{2}=0$ implies that $u(\cdot)\equiv 0$ and that $y(t)$ is in the null-space of $M(t)$, thus of $P(t)$, so $y(t_0)=0$ since $\scalidx{y(t_0)}{J_0 y(t_0)}{H}=0$.} Indeed, if $M(\cdot)\equiv 0$, and $J_0=0$, then \eqref{eq:2-3} is only a semi-norm, suited for trajectories with null initial condition, a subspace of $\cH$ which we will come back to in \Cref{sec:kernel_decomposition}.  For applications, $J_0$ can be taken very small and does not change the objective function if $y(t_0)$ is given. 
	
	\subsection{Identifying the kernel of LQ control}
	
	The main result is the following 
	\begin{thm}
		\label{theo3-1} We assume the coercivity of the drift \eqref{eq:2-1}, the strong convexity of the objective \eqref{eq:2-4}, and the invertibility \eqref{eq:invertibility} conditions. Define the family of operators $K(s,t)\in\cL(H,H)$ by the formula, with $\Phi_{A,P}(t,s)$ defined in \eqref{eq:def_PhiAP},
		\begin{align}
			&K(s,t):=\Phi_{A,P}(s,0)(J_0+P(t_0))^{-1}\Phi_{A,P}^{*}(t,0)\label{eq:3-3}\\
			&\hspace{2mm}+\int_{t_0}^{\min(s,t)}\Phi_{A,P}(s,\tau)B(\tau)N^{-1}(\tau)B^{*}(\tau)\Phi_{A,P}^{*}(t,\tau)d\tau.\nonumber
		\end{align}
		Then the space $(\cH,\noridx{\cdot}{\cH})$ defined by (\ref{eq:3-1}),(\ref{eq:3-2})
		is a RKHS associated with the kernel $K$ .
	\end{thm}
	
	\begin{proof}
		Fix $t\in[t_0,T]$, we set 
		\begin{equation}
			y_{zt}(s)=K(s,t)z\label{eq:3-4}
		\end{equation}
		and define the functions $\chi_{zt}(s)$ and $\rho_{zt}(s)$ for any $s\in[t_0,t]$ by the equation 
		\begin{gather}
			-\tfrac{d\chi_{zt}}{ds}+A^{*}(s)\chi_{zt}(s)=0,\nonumber\\ \chi_{zt}(t)=z,\label{eq:3-5} \\
			-\tfrac{d\rho_{zt}}{ds}+(A^{*}(s)+P(s)B(s)N^{-1}(s)B^{*}(s))\rho_{zt}(s)=0,\nonumber\\ \rho_{zt}(t)=z,\label{eq:3-6}
		\end{gather}
		which means also
		\begin{align}
			\chi_{zt}(s)&=\Phi_{A}^{*}(t,s)z\label{eq:3-4b},\\
			\rho_{zt}(s)&=\Phi_{A,P}^{*}(t,s)z\label{eq:3-7}.
		\end{align}
		From (\ref{eq:3-3}) we can then write, denoting by $\tau\in\R\mapsto\1_{\tau<t}$ the characteristic function of $(-\infty,t)$,
		\begin{align*}
			y_{zt}(s)&=\Phi_{A,P}(s,0)(J_0+P(t_0))^{-1}\rho_{zt}(t_0)\\
			&+\int_{t_0}^{s}\Phi_{A,P}(s,\tau)B(\tau)N^{-1}(\tau)B^{*}(\tau)\rho_{zt}(\tau)\1_{\tau<t}d\tau,
		\end{align*}
		and from the definition of $\Phi_{A,P}(s,\tau)$, we deduce that
		\begin{gather}
			\tfrac{dy_{zt}(s)}{ds}+(A(s)+B(s)N^{-1}(s)B^{*}(s)P(s))y_{zt}(s)\nonumber \\
			=B(s)N^{-1}(s)B^{*}(s)\rho_{zt}(s)\1_{s<t},\label{eq:3-8}\\
			y_{zt}(t_0)=(J_0+P(t_0))^{-1}\rho_{zt}(t_0). \nonumber
		\end{gather}
		Define next 
		\begin{equation}
			r_{zt}(s):=-(\rho_{zt}(s)-\chi_{zt}(s))\1_{s<t}\label{eq:3-9}.
		\end{equation}
		We note that $r_{zt}(\cdot)$ is continuous. Combining (\ref{eq:3-5}) and (\ref{eq:3-6}) we see easily that $r_{zt}(s)$ is the solution of the equation 
		\begin{gather}
			-\tfrac{dr_{zt}}{ds}+(A^{*}(s)+P(s)B(s)N^{-1}(s)B^{*}(s))r_{zt}(s)\nonumber\\
			=P(s)B(s)N^{-1}(s)B^{*}(s)\chi_{zt}(s)\1_{s<t},\label{eq:3-10}\\\textit{}
			\quad r_{zt}(T)=0.\nonumber
		\end{gather}
		We then introduce the function 
		\begin{equation}
			q_{zt}(s)=P(s)y_{zt}(s)+r_{zt}(s),\label{eq:3-11}
		\end{equation}
		whence 
		\begin{equation*}
			q_{zt}(t_0)=P(t_0)y_{zt}(t_0)+\chi_{zt}(t_0)-(J_0+P(t_0))y_{zt}(t_0).%
		\end{equation*}
		Using the relation (\ref{eq:2-80}) as well as (\ref{eq:3-8}),(\ref{eq:3-10}),
		we check easily that the pair $y_{zt}(s),q_{zt}(s)$ satisfies the
		system of forward backward differential equations 
		\begin{gather}
			\tfrac{dy_{zt}(s)}{ds}+A(s)y_{zt}(s)+B(s)N^{-1}(s)B^{*}(s)q_{zt}(s)\nonumber\\=B(s)N^{-1}(s)B^{*}(s)\chi_{zt}(s)\1_{s<t},\label{eq:3-12}\\
			-\tfrac{dq_{zt}(s)}{ds}+A^{*}(s)q_{zt}(s)=M(s)y_{zt}(s), \nonumber\\
			J_0y_{zt}(t_0)=\chi_{zt}(t_0)-q_{zt}(t_0) ,\quad q_{zt}(T)=0. \nonumber
		\end{gather}
		We can now proceed with the proof that $\cH$ is a RKHS associated
		to the kernel $K$ defined by (\ref{eq:3-3}). We have to check the properties \eqref{eq:3-100}-\eqref{eq:rep-prop}. The first property is clear by \eqref{eq:3-12}. Indeed $y_{zt}(\cdot)$ belongs to $\cH$ with control:
		\begin{equation*}
			u_{zt}(s)=-N^{-1}(s)B^{*}(s)(q_{zt}(s)-\chi_{zt}(s)\1_{s<t}),%
		\end{equation*}
		the latter satisfying \eqref{def_u_as_X-X'}. Indeed, by \eqref{eq:3-12}, $B(s)^{\ominus}[\tfrac{dy}{ds} + A(s) y(s)]=B(s)^{\ominus}B(s)N^{-1}(s)B^{*}(s)(\chi_{zt}(s)\1_{s<t}-q_{zt}(s))$. However $B(s)^{\ominus}B(s)$ is precisely the projector on $\operatorname{Ker}B(s)=\operatorname{Im}B^{*}(s)$ for the inner product $\scalidx{N(t)\cdot}{\cdot}{U}$. Consequently $B(s)^{\ominus}B(s)N^{-1}(s)B^{*}(s)=N^{-1}(s)B^{*}(s)$.
		
		It remains to prove \eqref{eq:rep-prop}. Take any
		$y(\cdot)\in\cH$, with $\cH$ defined in (\ref{eq:3-1}), then we have that
		\begin{align*}
			&\scalidx{y(\cdot)}{y_{zt}(\cdot)}{\cH}\\
			&\stackrel{\eqref{eq:3-2}}{=} \scalidx{y(t_0)}{J_0 y_{zt}(t_0)}{H}+\int_{t_0}^{T} \scalidx{M(s)y_{zt}(s)}{y(s)}{H}ds \nonumber\\
			&\hspace{2mm}-\int_{t_0}^{T}\scalidx{N(s)u(s)}{N^{-1}(s)B^{*}(s)(q_{zt}(s)-\chi_{zt}(s)\1_{s<t})}{U}ds\nonumber\\
			&=
			\scalidx{y(t_0)}{\chi_{zt}(t_0)-q_{zt}(t_0)}{H}+\int_{t_0}^{T} \scalidx{M(s)y_{zt}(s)}{y(s)}{H}ds\\
			&\hspace{5mm}-\int_{t_0}^{T} \scalidx{B(s)u(s)}{q_{zt}(s)-\chi_{zt}(s)\1_{s<t}}{H}ds\nonumber\\
			&\stackrel{\eqref{eq:3-12}}{=}
			\scalidx{y(t_0)}{\chi_{zt}(t_0)-q_{zt}(t_0)}{H}\\
			&\hspace{5mm}+\int_{t_0}^{T} \scalidx{-\tfrac{dq_{zt}(s)}{ds}+A^{*}(s)q_{zt}(s)}{y(s)}{H}ds\nonumber\\
			&\hspace{5mm}-\int_{t_0}^{T} \scalidx{\tfrac{dy(s)}{ds}+A(s)y(s)}{q_{zt}(s)-\chi_{zt}(s)\1_{s<t}}{H}ds.\nonumber
		\end{align*}
		Integrating by parts, recalling that $\chi_{zt}(s)=\Phi_{A}^{*}(t,s)z$ by \eqref{eq:3-4b}, we obtain by applying the variation of constants that
		\begin{align*}
			\scalidx{y(\cdot)}{y_{zt}(\cdot)}{\cH}&=\scalidx{y(t_0)}{\chi_{zt}(t_0)}{H}\nonumber\\ &\hspace{5mm}+\int_{t_0}^{T} \scalidx{B(s)u(s)}{\chi_{zt}(s)}{H}ds\\
			&=\scalidx{y(t)}{z}{H},
		\end{align*}
		which completes the proof of (\ref{eq:3-100}) and of the theorem.
	\end{proof}

	\subsection{Decomposition of the kernel into null-control and null-initial condition}\label{sec:kernel_decomposition}
	From now on, we denote $\cH$ by $\HK$. Equation \eqref{eq:3-3} induces us to split the kernel $K$ into 
	\begin{equation}
		K(s,t)=K^{0}(s,t)+K^{1}(s,t)\label{eq:3-13}
	\end{equation}
	with 
	\begin{align}
		K^{0}(s,t)&:=\Phi_{A,P}(s,0)(J_0+P(t_0))^{-1}\Phi_{A,P}^{*}(t,0),\label{eq:3-14}\\
		K^{1}(s,t)&:=\nonumber\\
		&\hspace{-5mm}\int_{t_0}^{\min(s,t)}\Phi_{A,P}(s,\tau)B(\tau)N^{-1}(\tau)B^{*}(\tau)\Phi_{A,P}^{*}(t,\tau)d\tau.\nonumber
	\end{align}
	We will see in \Cref{sec:solving_LQR} that $K^{1}$ is instrumental for the LQR. If we define $\xi_{zt}^{0}(s)=K^{0}(s,t)z$,  then we can check that
	the pair $\xi_{zt}^{0}(s)$$, \:\eta_{zt}^{0}(s)$ is the solution
	of the system 
	\begin{gather}
		\tfrac{d\xi_{zt}^{0}(s)}{ds}+A(s)\xi_{zt}^{0}(s)+B(s)N^{-1}(s)B^{*}(s)\:\eta_{zt}^{0}(s)=0,\nonumber\label{eq:3-15}\\
		-\tfrac{d\eta_{zt}^{0}(s)}{ds}+A^{*}(s)\eta_{zt}^{0}(s)=M(s)\xi_{zt}^{0}(s),\nonumber\\
		\xi_{zt}^{0}(t_0)=-\eta_{zt}^{0}(t_0)+\rho_{zt}(t_0). \nonumber
	\end{gather}
	with $\eta_{zt}^{0}(s)=P(s)\xi_{zt}^{0}(s)$. Similarly, setting $\xi_{zt}^{1}(s)=K^{1}(s,t)z$,
	we have the relations
	\begin{gather}
		\tfrac{d\xi_{zt}^{1}(s)}{ds}+A(s)\xi_{zt}^{1}(s)+B(s)N^{-1}(s)B^{*}(s)\:\eta_{zt}^{1}(s)\nonumber\\
		=B(s)N^{-1}(s)B^{*}(s)\chi_{zt}(s)\1_{s<t},\nonumber\label{eq:3-16}\\
		-\tfrac{d\eta_{zt}^{1}(s)}{ds}+A^{*}(s)\eta_{zt}^{1}(s)=M(s)\xi_{zt}^{1}(s),\nonumber\\
		\xi_{zt}^{1}(t_0)=0, \quad \eta_{zt}^{1}(T)=0. \nonumber
	\end{gather}
	and  $\eta_{zt}^{1}(s)=P(s)\xi_{zt}^{1}(s)+r_{zt}(s)$. We have also 
	\begin{align}
		y_{zt}(s)&=\xi_{zt}^{0}(s)+\xi_{zt}^{1}(s)\label{eq:3-18},\\
		q_{zt}(s)&=\eta_{zt}^{0}(s)+\eta_{zt}^{1}(s). \nonumber
	\end{align}
	Consider the Hilbert subspace of $\HK^{1}$ of functions with initial value equal to $0$, namely 
	\begin{multline}
		\HK^{1}=\{y(\cdot)\, | \,\tfrac{dy}{dt}+A(t)y(t)=B(t)u(t),y(t_0)=0,\\ \text{with}\;u(\cdot)\in L^{2}(t_0,T;U)\}.\label{eq:3-19}
	\end{multline}
	equipped with the same norm $\noridx{\cdot}{\HK}$ defined in \eqref{eq:3-2}. Then we claim 
	\begin{prop}
		\label{prop3-1} The Hilbert space \textup{$\HK^{1}$} is a RKHS associated with the operator-valued kernel $K^{1}(s,t).$
	\end{prop}
	
	\begin{proof}
		The proof is similar to the one of Theorem \ref{theo3-1} for $\cH$ and $K$. Therefore, we can use the notation $\HK^{1}=\cH_{K^{1}}$. Note that we can write 
		\begin{equation}
			\HK=\HK^{0}\oplus\HK^{1}\label{eq:3-21}
		\end{equation}
		where 
		\begin{equation}%
			\HK^{0}=\{y(\cdot)\, | \,\tfrac{dy}{dt}+A(t)y(t)=0\}.\label{eq:3-22}
		\end{equation}
	\end{proof}
	\noindent Conversely we do not have $\HK^{0}=\cH_{K^{0}}$, since for $y(\cdot)\in \HK^{0}$ 
	\begin{equation*}
		\scalidx{\xi_{zt}^{0}(\cdot)}{y(\cdot)}{\HK}\neq \scalidx{y(t)}{z}{H}.
	\end{equation*}
	This discrepancy between subspaces can be made more explicit by studying their projections.
	\subsection{Projections on subspaces}
	
	We have the property 
	\begin{prop}
		\label{prop3-2} The projection of $K(\cdot,t)z$ on \textup{$\cH_{K^{1}}$
			is $K^{1}(\cdot,t)z.$}
	\end{prop}
	
	\begin{proof}
		We need to check that 
		\begin{equation*}
			\scalidx{K(\cdot,t)z}{y(\cdot)}{\HK}=\scalidx{K^1(\cdot,t)z}{y(\cdot)}{\HK},\:\forall y(\cdot)\in\cH_{K^{1}}.
		\end{equation*}
		This is equivalent to proving that, for all $y(\cdot)\in\cH_{K^{1}}$, $\scalidx{\xi_{zt}^{0}(\cdot)}{y(\cdot)}{\HK}=\scalidx{K^0(\cdot,t)z}{y(\cdot)}{\HK}=0$. The latter is checked easily with techniques already used.
	\end{proof}
	We want also to find the projection of $K(\cdot,t)z$ on $\HK^{0}$.
	Let $\pi(t)\in\cL(H,H)$ be the solution of the linear equation 
	\begin{gather}
		-\tfrac{d\Psi}{dt}+A^{*}(t)\Psi(t)+\pi(t)(\tfrac{d\varphi}{dt}+A(t)\varphi(t))=M(t)\varphi(t)\label{eq:3-24}\\
		\text{for all }\varphi(\cdot)\text{ such that }\tfrac{d\varphi}{dt}+A(t)\varphi(t)\in L^{2}(t_0,T;H), \nonumber\\
		\quad\Psi(t)=\pi(t)\varphi(t), \; \pi(T)=0.\nonumber
	\end{gather}
	Note that if $M(\cdot)\equiv 0$, then $P(\cdot)\equiv 0$ and $\pi(\cdot)\equiv 0$. We have the following result
	\begin{prop}
		\label{prop3-3} The projection of $K(\cdot,t)z$ on $\HK^{0}$ is $\Phi_{A}(\cdot,0)(J_0+\pi(t_0)$)$^{-1}\Phi_{A}^{*}(t,0)z$.
	\end{prop}
	\begin{proof} Recall that, by definition \eqref{eq:3-4b}, $\chi_{zt}(s)=\Phi_{A}^{*}(t,s)z$ for $s<t$. We claim that
		the projection of $K(\cdot,t)z$ is $\Phi_{A}(\cdot,0)(J_0+\pi(t_0))^{-1}\chi_{zt}(t_0).$
		Let $\zeta(\cdot)$ be the projection of $K(\cdot,t)z$ on $\HK^{0}.$
		We have 
		\begin{equation*}
			\tfrac{d\zeta}{ds}+A(s)\zeta(s)=0 ,\quad \zeta(t_0)=\zeta_{0},
		\end{equation*}
		and the problem is to find $\zeta_{0}$. The projection has to satisfy the optimality condition 
		\begin{equation*}
			\scalidx{K(\cdot,t)z}{y(\cdot)}{\HK}=\scalidx{\zeta(\cdot)}{y(\cdot)}{\HK},\;\forall y(\cdot)\in\HK^{0}.
		\end{equation*}
		Since the trajectories of $\HK^{0}$ have null control, we obtain after an integration by parts, since $q_{zt}(T)=0$, that
		\begin{align*}
			&\scalidx{K(\cdot,t)z}{y(\cdot)}{\HK}\nonumber\\
			&=\scalidx{\chi_{zt}(t_0)-q_{zt}(t_0)}{y(t_0)}{H}\nonumber\\&\hspace{1cm}+\int_{t_0}^{T}\scalidx{M(s)y_{zt}(s)}{y(s)}{H}ds \nonumber\\
			&=\scalidx{\chi_{zt}(t_0)-q_{zt}(t_0)}{y(t_0)}{H}\nonumber\\&\hspace{1cm}+\int_{t_0}^{T}\scalidx{-\tfrac{dq_{zt}(s)}{ds}+A^{*}(s)q_{zt}(s)}{y(s)}{H}ds\nonumber\\
			&=\scalidx{\chi_{zt}(t_0)-q_{zt}(t_0)}{y(t_0)}{H}\nonumber\\&\hspace{5mm}+\int_{t_0}^{T}\scalidx{-\tfrac{dq_{zt}(s)}{ds}}{y(s)}{H}+\scalidx{q_{zt}(s)}{-\tfrac{dy(s)}{ds}}{H}ds\nonumber\\
			&=\scalidx{\chi_{zt}(t_0)}{y(t_0)}{H}.\nonumber
		\end{align*}
		On the other hand, we have also 
		\begin{align*}
			\scalidx{\zeta(\cdot)}{y(\cdot)}{\HK}&=\scalidx{J_0\zeta(t_0)}{y(t_0)}{H}\nonumber\\ &\hspace{5mm}+\int_{t_0}^{T} \scalidx{M(s)\zeta(s)}{y(s)}{H}ds.
		\end{align*}
		If we introduce $\eta(s)$, solution of 
		\[
		-\tfrac{d\eta}{ds}+A^{*}(s)\eta(s)=M(s)\zeta(s), \quad \eta(T)=0,
		\]
		we obtain the relation 
		\[
		\scalidx{\zeta(\cdot)}{y(\cdot)}{\HK}=\scalidx{J_0\zeta(t_0)}{y(t_0)}{H}+\scalidx{\eta(t_0)}{y(t_0)}{H}.
		\]
		The optimality condition reduces to 
		\[
		\scalidx{\chi_{zt}(t_0)}{y_0}{H}=\scalidx{J_0\zeta(t_0)}{y_0}{H}+\scalidx{\eta(t_0)}{y_0}{H},\forall y_{0}\in H.
		\]
		Therefore we obtain 
		\[
		J_0\zeta_{0}=\chi_{zt}(t_0)-\eta(t_0).
		\]
		However, we can check easily that $\eta(s)=\pi(s)\zeta(s).$ Therefore
		$\zeta_{0}=(J_0+\pi(t_0))^{-1}\chi_{zt}(t_0)$ which yields the result.
	\end{proof}

	\subsection{Example of heat equation with distributed control}\label{sec:heat_example}
	As discussed above, we here focus on bounded $B(\cdot)\in L^\infty$ and parabolic equations (unbounded/hyperbolic would require a few changes). Take $V=H^1(\R^d,\R)$, $H=L^2(\R^d,\R)$, $A(\cdot)\equiv-\Delta$ and $B(\cdot)\equiv\Id_H$, then the heat equation with distributed control writes as
	\begin{equation}
		\tfrac{dy}{dt}=\Delta y(t)+u(t), \, y(t_0)=y_{0}\in H\label{eq:heat}.
	\end{equation}
	This equation can sketchily describe a microwave oven (local heating) in a refrigerator (global cooling), making the operators time-dependent would account for a turntable. As objective, let us take $J_0=\lambda\Id_H$ with $\lambda>0$, $M(\cdot)\equiv 0$ and $N(\cdot)\equiv \Id_H$, thus $P(\cdot)\equiv 0$, and $\Phi_{A,P}(t,s)=\Phi_{A}(t,s)$. In this well-known context, the (integral) operator $\Phi_{A}(t,s)=e^{-A(t-s)}$ is merely the heat semi-group associated to the heat kernel, for $t>s$,
	 $$k(t-s,x,y)=\dfrac{1}{(4\pi(t-s))^{\tfrac{d}{2}}}e^{-\tfrac{\|x-y\|^2_d}{4(t-s)}}.$$
	Using that $A$ is self-adjoint and the known expression of the Fourier transform of a normalized Gaussian, one can show that $\int_0^{2s}k(\tau,x,y) d\tau= k(s^2,x,y)$ and consequently that,  for $t>s$, $K^1(s,t)=\frac{1}{2}[\int_0^{2s}e^{-A\tau} d\tau]\circ e^{-A(t-s)}$ is a kernel integral operator with kernel $k(t-s+s^2,x,y)/2$. On the other hand $K^0(s,t)=e^{-A(t+s)}/\lambda$ has for kernel $k(t+s,x,y)/\lambda$. This allows for explicit handling of the kernel $K$ in applied cases with various objective functions, as considered in the next section.
	
	\section{Solving control problems}\label{sec:optimal_control}
	In this section, we return to optimal control problems, first with only a differentiable terminal cost, of which the LQR is a special case, and then for more general objective functions.
	\subsection{Final nonlinear term - Mayer problem}
	We consider the dynamic system 
	\begin{equation}
		\tfrac{dy}{dt}+A(t)y(t)=B(t)u(t) , \quad y(t_0)=y_{0}.\label{eq:4-1}
	\end{equation}
	We want to find the pair $y_{0},u(\cdot)$ in order to minimize 
	\begin{multline}
		J(u(\cdot),y_{0}):=g(y(T))+\tfrac{1}{2}\scalidx{y(t_0)}{J_0 y(t_0)}{H}\\
		\hspace{5mm}+\tfrac{1}{2}\int_{t_0}^{T} \big[\scalidx{M(t)y(t)}{y(t)}{H}+ \scalidx{N(t)u(t)}{u(t)}{U}\big]dt,\label{eq:4-2}
	\end{multline}
	where $h\mapsto g(h)$ is a Gâteaux differentiable function on
	$H.$ Using the norm $\noridx{\cdot}{\HK}$ defined in \eqref{eq:3-2}, this problem can be formulated as minimizing a functional on $\HK$, namely 
	\begin{equation}
		\mathcal{J}(y(\cdot)):=g(y(T))+\tfrac{1}{2}\noridx{y(\cdot)}{\HK}^{2}.\label{eq:4-3}
	\end{equation}
	If $\hat{y}(\cdot)$ is a minimizer, it satisfies the Euler equation 
	\begin{equation}
		\scalidx{Dg(\hat{y}(T))}{\zeta(T)}{H}+\scalidx{\hat{y}(\cdot)}{\zeta(\cdot)}{\HK}=0,\:\forall\zeta(\cdot)\in\HK.\label{eq:4-4}
	\end{equation}
	But we may write by the reproducing property \eqref{eq:rep-prop}
	\[
	\scalidx{Dg(\hat{y}(T))}{\zeta(T)}{H}=\scalidx{K(\cdot,T)Dg(\hat{y}(T)}{\zeta(\cdot)}{\HK}
	\]
	and (\ref{eq:4-4}) yields immediately the equation for $\hat{y}(\cdot)$
	\begin{equation}
		K(\cdot,T)Dg(\hat{y}(T))+\hat{y}(\cdot)=0\label{eq:4-5}.
	\end{equation}

	\subsection{Recovering the standard solution of the LQR}\label{sec:solving_LQR}
	We can now go back to the standard LQR problem \eqref{eq:2-2}-\eqref{eq:2-3}, where the initial state $y_{0}$ is known. The state $y(\cdot)$ can
	be written as follows 
	
	\begin{equation}
		y(s)=\Phi_{A}(s,0)y_{0}+\zeta(s)\label{eq:4-6}
	\end{equation}
	where $\zeta(\cdot)$ satisfies 
	
	\begin{equation*}
		\tfrac{d\zeta}{ds}+A(s)\zeta(s)=B(s)u(s) , \quad \zeta(t_0)=0.\label{eq:4-7}
	\end{equation*}
	Therefore $\zeta(\cdot)\in\cH_{K^{1}}.$ We write $y_{0}(s)=\Phi_{A}(s,0)y_{0}$.
	The cost (\ref{eq:2-3}) becomes 
	\begin{align*}
		J(u(\cdot))&=\int_{t_0}^{T} \scalidx{M(t)y_0(t)}{y_0(t)}{H}dt\\
		&+\int_{t_0}^{T} \scalidx{M(t)\zeta(t)}{\zeta(t)}{H}dt\\
		&\hspace{-1cm}+2\int_{t_0}^{T} \scalidx{M(t)y_0(t)}{\zeta(t)}{H}dt+\int_{t_0}^{T} \scalidx{N(t)u(t)}{u(t)}{U}dt.%
	\end{align*}
	The problem amounts to minimizing 
	\begin{equation}
		\mathcal{J}(\zeta(\cdot))=\noridx{\zeta(\cdot)}{\HK}^{2}+2\int_{t_0}^{T} \scalidx{M(t)y_0(t)}{\zeta(t)}{H}dt,\label{eq:4-9}
	\end{equation}
	on the Hilbert space $\cH_{K^{1}}.$ The functional $\mathcal{J}(\zeta(\cdot))$ can also be written as 
	\begin{equation}
		\mathcal{J}(\zeta(\cdot))=\noridx{\zeta(\cdot)}{\HK}^{2}+2 \scalidx{\zeta(\cdot)}{\int_{t_0}^{T}K^{1}(\cdot,t)M(t)y_{0}(t)dt}{H},\label{eq:4-10}
	\end{equation}
	and the minimizer is obtained immediately by the formula 
	\begin{equation}
		\hat{\zeta}(s)=-\int_{t_0}^{T}K^{1}(s,t)M(t)y_{0}(t)dt.\label{eq:4-11}
	\end{equation}
	Using the second formula (\ref{eq:3-14}) we write 
	\begin{align*}
		\hat{\zeta}(s)&=-\int_{t_0}^{T}\big[\int_{t_0}^{\min(s,t)}\Phi_{A,P}(s,\tau)B(\tau)N^{-1}(\tau)\\
		&\hspace{2cm}\times B^{*}(\tau)\Phi_{A,P}^{*}(t,\tau)d\tau\big]M(t)y_{0}(t)dt\\
		&=-\int_{t_0}^{s}\Phi_{A,P}(s,\tau)B(\tau)N^{-1}(\tau)B^{*}(\tau)\\
		&\hspace{2cm}\times\big(\int_{\tau}^{T}\Phi_{A,P}^{*}(t,\tau)M(t)y_{0}(t)dt\big)d\tau.
	\end{align*}
	One checks easily that 
	\[
	\int_{\tau}^{T}\Phi_{A,P}^{*}(t,\tau)M(t)y_{0}(t)dt=P(\tau)y_{0}(\tau).
	\]
	Therefore $\hat{\zeta}(s)$ satisfies 
	\begin{gather*}
		\tfrac{d\hat{\zeta}(s)}{ds}+(A(s)+B(s)N^{-1}(s)B^{*}(s)P(s))\hat{\zeta}(s)\\
		=-B(s)N^{-1}(s)B^{*}(s)P(s)y_{0}(s).
	\end{gather*}
	and $y(s)=y_{0}(s)+\hat{\zeta}(s)$ satisfies 
	\begin{equation}
		\tfrac{dy}{ds}+(A(s)+B(s)N^{-1}(s)B^{*}(s)P(s))y(s)=0 , \quad y(t_0)=0,\label{eq:4-12}
	\end{equation}
	which coincides with the evolution equation \eqref{eq:2-9}-\eqref{eq:2-10} of the optimal state
	for problem \eqref{eq:2-2}-\eqref{eq:2-3}.

	\subsection{More general objectives: state constraints and intermediary points}
	
	More generally, within the same framework, for any collection of time points $(t_n)_{n=1}^N\in[t_0,T]^N$, we can tackle objectives of the form
	\begin{equation}
		\mathcal{J}(y(\cdot))=L((y(t_n))_{n=1}^N,\noridx{y(\cdot)}{\HK}^{2})\label{eq:gen_objective}
	\end{equation}%
	for a given extended-valued function $L:H^N\times[0,+\infty]\rightarrow\R\cup\{+\infty\}$. Notably the LQR problem \eqref{eq:2-3} with terminal cost is a special case of \eqref{eq:gen_objective}. Moreover \eqref{eq:gen_objective} can incorporate indicator functions of constraints sets over $y(\cdot)$, and thus handle a finite number state constraints, beside a terminal constraint as considered in \cite{Przyuski2014}. For an infinite number of state constraints, e.g.\ holding on the whole time interval $[t_0,T]$, one can extend the technique of \cite{aubin2020hard_control}. Indeed, owing to our rewriting of optimal control problems with operator-valued kernels, we have access to ``representer theorems`` such as:
	\begin{thm}[Theorem 4.2, \cite{Micchelli2005}]\label{thm:rep}
		If for every $\b z \in H^N$ the function $h: \xi\in\R_{+} \mapsto L(\b z, \xi)\in\R_{+}\cup\{+\infty\}$ is strictly increasing and $\hat y(\cdot) \in \HK$ minimizes the functional \eqref{eq:gen_objective}, then $\hat y(\cdot)=\sum_{n=1}^N K(\cdot,t_n) z_{n}$ for some $\left\{z_{n}\right\}_{n=1}^N \subseteq H$. In addition, if $L$ is strictly convex, the minimizer is unique.
	\end{thm}
	\Cref{thm:rep} can be seen as a generalization of \eqref{eq:4-4}, whereas \eqref{eq:4-5} characterizes the optimal $z_{T}=-Dg(\hat{y}(T))$. The optimal control can be obtained by linearity as $\hat u(\cdot)=\sum_{n=1}^N u_{z_{n}t_n}(\cdot)=:\sum_{n=1}^N U(\cdot,t_n)z_{n}$. We refer to \cite{aubin2020hard_control} for more on the matrix $U(\cdot,t)z=u_{zt}(\cdot)$.
	
	\section*{Conclusion}
	Through \Cref{theo3-1}, we have improved upon \cite{aubin2020hard_control}, proposing a more concise formula for the LQ kernel and extending the setting to the control of infinite dimensional linear systems. Curiously, this type of connection between control theory and kernel methods is quite new and we think that a lot can be done in this direction of research, especially considered the growing interest for the Koopman operator, see \cite{brunton2021modern} for a recent review, or for Model Predictive Control, see e.g.\ \cite{grune2020efficient}. The kernels we considered appear in linear-quadratic optimal control because of Hilbertian vector spaces of trajectories, while,  for estimation problems, the kernels appear through covariances of Gaussian processes as investigated in \cite{aubin2022Kalman}. It is this ``dual'', deterministic and stochastic, nature of kernels which underlies the ``duality'' between optimal control and estimation in the Linear-Quadratic case.
	
	\bibliographystyle{unsrt}
	\bibliography{oRKHS_InfDimControl}	

\begin{thebibliography}{10}

\bibitem{aubin2020hard_control}
Pierre-Cyril Aubin-Frankowski.
\newblock Linearly constrained linear quadratic regulator from the viewpoint of
  kernel methods.
\newblock {\em {SIAM} Journal on Control and Optimization}, 59(4):2693--2716,
  2021.

\bibitem{Aronszajn1943}
Nachman Aronszajn.
\newblock La th{\'{e}}orie des noyaux reproduisants et ses applications:
  {Premi{\`{e}}re Partie}.
\newblock {\em Mathematical Proceedings of the Cambridge Philosophical
  Society}, 39(3):133--153, 1943.

\bibitem{aronszajn50theory}
Nachman Aronszajn.
\newblock Theory of reproducing kernels.
\newblock {\em Transactions of the American Mathematical Society}, 68:337--404,
  1950.

\bibitem{berlinet04reproducing}
Alain Berlinet and Christine Thomas-Agnan.
\newblock {\em Reproducing Kernel Hilbert Spaces in Probability and
  Statistics}.
\newblock Kluwer, 2004.

\bibitem{steinwart08support}
Ingo Steinwart and Andreas Christmann.
\newblock {\em Support Vector Machines}.
\newblock Springer, 2008.

\bibitem{aubin2022Kalman}
Pierre-Cyril Aubin-Frankowski and Alain Bensoussan.
\newblock The reproducing kernel {Hilbert} spaces underlying linear {SDE}
  estimation, {Kalman} filtering and their relation to optimal control.
\newblock {\em Pure and Applied Functional Analysis}, 2022.
\newblock (to appear) \url{https://arxiv.org/abs/2208.07030}.

\bibitem{Micchelli2005}
Charles~A. Micchelli and Massimiliano Pontil.
\newblock On learning vector-valued functions.
\newblock {\em Neural Computation}, 17(1):177--204, January 2005.

\bibitem{caponetto2008universal}
Andrea Caponnetto, Charles~A. Micchelli, Massimiliano Pontil, and Yiming Ying.
\newblock Universal multi-task kernels.
\newblock {\em Journal of Machine Learning Research}, 9(52):1615--1646, 2008.

\bibitem{kadri2016}
Hachem Kadri, Emmanuel Duflos, Philippe Preux, St{{\'e}}phane Canu, Alain
  Rakotomamonjy, and Julien Audiffren.
\newblock Operator-valued kernels for learning from functional response data.
\newblock {\em Journal of Machine Learning Research}, 17(20):1--54, 2016.

\bibitem{lions1971optimal}
Jacques-Louis Lions.
\newblock {\em Optimal control of systems governed by partial differential
  equations}.
\newblock Springer-Verlag, Berlin,New York, 1971.

\bibitem{lions1988controlabilite}
Jacques Lions.
\newblock {\em Contr\^{o}labilit\'{e} exacte, perturbations et stabilisation de
  systemes distribu\'{e}s}.
\newblock Masson, Paris, 1988.

\bibitem{Bensoussan2007}
Alain Bensoussan, Giuseppe~Da Prato, Michel~C. Delfour, and Sanjoy~K. Mitter.
\newblock {\em Representation and Control of Infinite Dimensional Systems}.
\newblock Birkh\"{a}user Boston, 2007.

\bibitem{Pritchard1987}
A.~J. Pritchard and D.~Salamon.
\newblock The linear quadratic control problem for infinite dimensional systems
  with unbounded input and output operators.
\newblock {\em {SIAM} Journal on Control and Optimization}, 25(1):121--144,
  January 1987.

\bibitem{Senkene1973}
{\'{E}}.~Senkene and A.~Tempel{\textquotesingle}man.
\newblock Hilbert spaces of operator-valued functions.
\newblock {\em Mathematical Transactions of the Academy of Sciences of the
  Lithuanian {SSR}}, 13(4):665--670, October 1973.

\bibitem{burbea1984banach}
Jacob Burbea.
\newblock {\em Banach and Hilbert spaces of vector-valued functions : their
  general theory and applications to holomorphy}.
\newblock Pitman, Boston, 1984.

\bibitem{Przyuski2014}
K.Maciej Przy{\l}uski.
\newblock On an infinite dimensional linear-quadratic problem with fixed
  endpoints: The continuity question.
\newblock {\em International Journal of Applied Mathematics and Computer
  Science}, 24(4):723--733, December 2014.

\bibitem{brunton2021modern}
Steven~L. Brunton, Marko Budisic, Eurika Kaiser, and J.~Nathan Kutz.
\newblock Modern {Koopman} theory for dynamical systems, 2021.

\bibitem{grune2020efficient}
Lars Grune, Manuel Schaller, and Anton Schiela.
\newblock Efficient {MPC} for parabolic {PDEs} with goal oriented error
  estimation.
\newblock {\em SIAM Journal on Scientific Computing (to appear)}, 2021.

\end{thebibliography}
\end{document}